\documentclass[12pt]{article} 
\usepackage {amsmath, amsfonts, amssymb, latexsym}
\usepackage{xy} 
\xyoption{all} 

\newtheorem{theorem}{Theorem}
\newtheorem{definition}{Definition}
\newtheorem{lemma}[theorem]{Lemma}
\newtheorem{proposition}[theorem]{Proposition}
\newtheorem{example}{Example}

\def\U{\mathcal{U}}
\def\B{\mathcal{B}}
\def\R{\mathbb{R}}
\def\su{\mathfrak{su}}
\def\sl{\mathfrak{sl}}

\def\calO{\mathcal{O}}
\def\d{\partial} 
\def\tri{\triangle} 
\def\id{\operatorname{id}}
\def\ad{\operatorname{ad}}
\def\Mapsto{\longmapsto} 
\def\A{\mathfrak{A}}
\def\H{\mathcal{H}}
\def\C{\mathbf{C}}
\def\g{\mathfrak{g}}
\def\h{\mathfrak{h}}
\def\o{\mathop\otimes} 
\def\t{\mathbf{t}}

\def\qbinom#1#2#3{\left[ \begin{matrix} #1 \\ #2 \end{matrix} \right]_{#3}}
\def\Cal#1{\mathcal{#1}}

\def\To{\longrightarrow}

\def\eps{\epsilon} 
\def\Z{{\mathbb Z}}
\def\abs#1{\left| #1 \right|}

\def\prooflem#1{\vskip 0.05 in 
        \noindent {\it Proof of Lemma \ref{#1}.} $\ \ $ } 
 
\def\norm#1{\left\| #1 \right\|}

\begin{document}
\title{Introduction to Quantum Group Theory}
\date{May 2001} 
\author{ William Gordon Ritter \\ Jefferson Physical Laboratory \\ 
Harvard University, Cambridge, MA}
\maketitle 

\begin{abstract} 
This is a short, self-contained expository survey, 
focused on algebraic and analytic 
aspects of quantum groups. Topics covered include the definition of 
``quantum group,'' the Yang-Baxter equation, quantized universal enveloping 
algebras, representations of braid groups, the KZ equations and the 
Kohno-Drinfeld theorem, and finally compact quantum groups and 
the analogue of Haar measure for compact quantum groups. 
\end{abstract}

\section{Introduction} 
There are two main approaches to quantum group theory, the purely algebraic 
approach, and the $C^*$-algebra approach, which uses deep connections with 
functional analysis. A $C^*$ algebra is a Banach algebra satisfying 
$\norm{a^* \, a} = \norm{a}^2$ for all elements $a$. Equivalently, 
a $C^*$-algebra is a norm closed $\ast$-subalgebra of $\B(\H)$, where $\H$ denotes
a Hilbert space. 
The simple example of quantum matrix groups fits 
naturally into both frameworks.  In this paper we endeavor to 
\begin{enumerate} 
\item[$\bullet$] Define Hopf algebras, and their quasitriangular structures, which 
is the starting point of the algebraic approach. 
\item[$\bullet$] Discuss the Yang-Baxter equation and quantum matrix algebras, 
including an algebro-geometric view of the moduli of solutions. 
\item[$\bullet$] Give the construction of the quantized universal enveloping 
algebra in detail, both for the simplest nontrivial example ($\sl_2$) 
and the general finite-dimensional Lie algebra $\g$. 
\item[$\bullet$] Discuss the connection of quantum groups with representations 
of braid groups. 
\item[$\bullet$] Discuss an application of these techniques to 
understanding recent results of Drinfeld and Kohno regarding 
the Knizhnik-Zamolodchikov Equations. 
\item[$\bullet$]  Discuss compact quantum groups, Haar measure on these groups, 
and their representation theory. This can be viewed in light of Connes' 
formulation of noncommutative geometry. 
\end{enumerate}

\section{Hopf Algebras and the Universal R-matrix} 

\subsection{Notation} 
We begin by introducing the following \emph{leg numbering notation}, 
which is due to Sweedler. If $A$ is any unital 
algebra, and $R \in A \otimes A$, 
then in general $R$ admits a representation in the form 
$R = \sum \alpha_i \otimes \beta_i$ (finite sum) for 
$\alpha_i, \beta_i$ elements of $A$. There are exactly three ways of 
extending $R$ to an element of $A \otimes A \otimes A$ by tensoring 
with $1 \in A$. The leg numbering notation is a way of distinguishing 
between them: 
\begin{eqnarray*} 
R_{12} &:=& \sum_i \alpha_i \otimes \beta_i \otimes 1 \\ 
R_{13} &:=& \sum_i \alpha_i \otimes 1 \otimes \beta_i \\
R_{23} &:=& \sum_i 1 \otimes \alpha_i \otimes \beta_i \\
\end{eqnarray*} 
This notation is useful, in particular, for specifying the defining 
relations of the universal R-matrix, which we will do. 
More generally, we will write $R_{ij}$ for the extension of $R$ 
to $A^{\otimes n}$ which ``is'' $R$ in the $i$th and the $j$th component. 

In this paper, $\tau$ will always denote the flip operator 
which acts linearly on the second tensor power of a module 
by $\tau(a \otimes b) = b \otimes a$.

\subsection{Hopf Algebras}

We first define coalgebras and bialgebras. 
The starting point is to write $\mu : A \otimes A \to A$ 
for the multiplication map of an asspociative algebra over a field $k$, 
and to note that a unit for the algebra $A$ certainly 
defines a map $k \to A$ which we can call $\eta$. Then, 
the associativity of the multiplication and the nice multiplicative 
property of the unit can be written as follows: 
\begin{equation} \label{def:algebra} 
\begin{matrix} 
\text{associativity}  &  \text{unit}  \\ 
\xymatrix@C=7pt{ & A \o A \o A \ar[dr]^{\id \o \mu} \ar[dl]_{\mu \o \id}  & \\ 
A \o A \ar[dr]_{\mu} & & A \o A \ar[dl]^{\mu}  \\ 
& A & }  
& 
\xymatrix@C=7pt{ A \o A \ar[dr]^{\mu}  & \\  k \otimes A \ar[u]^{\eta \o \id} \ar@{=}[r] & A } 
\quad 
\xymatrix@C=7pt{ A \o A \ar[dr]^{\mu}  & \\  A \otimes k \ar[u]^{\id \o \eta} \ar@{=}[r] & A } 
\end{matrix} 
\end{equation} 
We thus take the definition of an \emph{algebra} to be 
a vector space with the additional structure of a pair of 
linear maps $(\mu, \eta)$ making \eqref{def:algebra} 
a commutative diagram. 

We can now dualize this construction, so that a \emph{coalgebra} 
is a vector space $C$ with the additional structure of a pair of 
linear maps $(\tri, \eps)$ making \eqref{def:coalgebra} 
a commutative diagram: 
\begin{equation} \label{def:coalgebra} 
\begin{matrix} 
\text{coassociativity}  &  \text{counit}  \\ 
\xymatrix@C=7pt{ & C \o C \o C   & \\ 
C \o C \ar[ur]^{\tri\o\id}  & & C \o C \ar[ul]_{\id\o\tri}   \\ 
& C \ar[ur]_{\tri} \ar[ul]^{\tri} & }  
& 
\xymatrix@C=7pt{ C \o C \ar[d]_{\eps\o\id}   & \\  k \otimes C \ar@{=}[r]  & C \ar[ul]_{\tri} } 
\quad 
\xymatrix@C=7pt{ C \o C \ar[d]_{\id\o\eps}   & \\  C \otimes k \ar@{=}[r]  & C \ar[ul]_{\tri} } 
\end{matrix} 
\end{equation}

A bialgebra structure on a vector space $A = C$ 
is a quadruple of objects ($\mu, \eta, \tri, \eps$) 
which satisfy all of the commutative diagrams 
\eqref{def:algebra}-\eqref{def:coalgebra} as well as 
the following compatibility equations: 
\[
\tri(hg) = \tri(h) \tri(g), \quad \tri(1) = 1 \o 1, \quad 
\eps(hg) = \eps(h) \eps(g), \quad \eps(1) = 1, \qquad 
\text{ for all } g, h \in A 
\]
We define the \emph{convolution product} of two 
linear maps $f, g : A \To A$ by the formula 
$f \ast g = \mu (f \o g) \tri$. Not surprisingly, 
$\ast$ is associative and has neutral element $\eta \circ \eps$. 
In this situation, 
a bialgebra $(A, \mu, \eta, \tri, \eps)$ is said 
to be a \emph{Hopf Algebra} if the identity map 
$\id : A \to A$ is invertible for the convolution product, 
and its inverse is the \emph{antipode} $S$. The defining property 
of this antipode can also be represented as a commutative
diagram: 
\[
\xymatrix{ 
A \ar[r]^{\eps} \ar[d]^{\tri}  &  k \ar[r]^{\eta} & A \\ 
A \o A \ar[rr]^{\id \o S, \ \ S \o \id}  &  & A \o A \ar[u]_{\mu} }
\]

\begin{definition} \label{def:QT-hopf} 
A \emph{quasitriangular Hopf algebra}, also called a 
\emph{braided Hopf algebra}, consists of the data of a Hopf algebra, 
$(A,\mu,\eta,\tri,\eps,S)$ together with an invertible element 
$R \in A \otimes A$ satisfying the following two conditions: 
\begin{enumerate} 
\item[(i)]  $R \tri(x) R^{-1} = \tri'(x)$ for all $x \in A$. 

\item[(ii)]  $(\tri \otimes \id)(R) = R_{13} R_{23}$ and 
$(\id\otimes\tri)(R) = R_{13}R_{12}$. 
\end{enumerate} 
where $\tri' = \tau \circ \tri$ is the opposite comultiplication. 
In this situation, $R$ is called the \emph{universal R-matrix}. 
\end{definition}

\section{The Yang-Baxter Equation} 

The Yang-Baxter equation \eqref{eqn:YBE} is simple to 
write down, but it is more complicated
to understand the motivation for its study; hence, 
this deeper kind of understanding will form the 
subject matter of this section. 

\begin{lemma} \label{lemma:YBE} 
From axioms (i)-(ii) in Definition \ref{def:QT-hopf} it follows 
that $R$ satisfies the \emph{Yang-Baxter Equation}: 
\begin{equation} \label{eqn:YBE} 
R_{12} R_{13} R_{23} = R_{23} R_{13} R_{12} 
\end{equation} 
\end{lemma} 

\prooflem{lemma:YBE} 
Compute $(\id \o \tau \circ \tri) R$ in two ways, 
using the second part of (ii), or by first using (i) and 
then the second part of (ii). $\Box$. 

\begin{example} 
If $G$ is an algebraic group or monoid then we define 
$\calO(G)$ to be the bialgebra of polynomial functions 
on $G$. Of particular interest is $G = M(n) = n \times n $ 
matrices over a field. There is a ``standard'' quantization 
of $\calO(M(n))$ for which the quantizing Yang-Baxter matrix
is given by 
\[
R = \exp(t \gamma_q) \exp(t\beta) \exp(t\gamma_q) 
\]
where $\gamma_q = \sum_{i<j} e_{ij} \wedge e_{ji}$, and
$\beta$ is a generic element of $M(n) \otimes M(n)$. 
\end{example}

\subsection{Quantum Matrices and the Yang-Baxter Variety} 

Specializing to the matrix algebra $M_n = M_n(k)$, 
we fix an element $R \in M_n \otimes M_n$ satisfying \eqref{eqn:YBE}, 
(an ``R-matrix'' from now on) and we first define the bialgebra 
$A(R)$ to be generated by unit element 1, together with 
$n^2$ indeterminates $\t = \{ {t^i}_j \}$ with relations 
$R \t_1 \t_2 = \t_2 \t_1 R$ and coalgebra structure: 
\begin{equation} \label{eqn:A(R)} 
\tri \t = \t \otimes \t, \quad 
\eps \t = \id 
\end{equation} 
We resist the temptation to write everything using index notation 
on the tensors, but as an example, note that the first of equations 
\eqref{eqn:A(R)} is shorthand for $\tri {t^i}_j = {t^i}_a \otimes {t^a}_j$, 
in which the sum over $a$ is understood according to the summation convention. 
The bialgebra $A(R)$ is called the algebra of 
\emph{quantum matrices}.

It is interesting to consider the moduli space of all 
invertible elements 
$R \in M_n \otimes M_n$ satisfying \eqref{eqn:YBE}. In this situation, 
it makes sense to quotient by the equivalence of different normalizations. 
We define $YB_n$ to be the resulting moduli space. 
It has naturally the topology of an algebraic variety, since 
an R-matrix can be interpreted as an element of the common zero 
loci of a collection of polynomials over $k$. 
In this interpretation, the construction of the quantum 
matrix algebra $A(R)$ which was discussed just above forms a 
kind of bundle over this variety, but the variety $YB_n$ 
itself is far from being a smooth manifold. In general it contains 
singular disconnected points, lines, planes, and jumps in dimension
(even when $n$ is held fixed)!

\subsection{Quantization of Poisson brackets} 

Drinfeld studied the quantization of Poisson brackets on 
the commutative algebra $C^\infty(G)$ where $G$ is a Lie group. 
Let $X_i$ denote a basis of $\g = $ the Lie algebra of $G$. 
Let $\d_i$ be the right-invariant vector field on $G$ corresponding 
to $X_i$, and let $\d_i'$ be the corresponding left-invariant vector
field. Let $R \in \bigwedge^2(\g)$ be given by $R = r^{ij} X_i \otimes X_j$. 
The Poisson brackets originally studied by Sklyanin, and more carefully
by Drinfeld, are then given by equations of the form 
\[
\{f,g\} = r^{ij} (\d_i f \, \d_j g - \d_i' f\, \d_j' g), \qquad
f, g \in C^\infty(G) 
\]

\section{$\su_2$ and $\sl_2$ Quantum Groups} 
\label{sec:sl2} 

In a course on finite-dimensional complex Lie algebras, it often 
makes sense to first study detailed properties of the Lie algebra 
$\sl_2$, or its real form $\su_2$, because the root space decomposition
shows that a semisimple Lie algebra $\g$ can be viewed as containing
a number of copies of $\sl_2$ as subalgebras, and the properties 
of these subalgebras are crucial in deriving important properties of 
$\g$ and its representations. It turns out that something similar 
is true for the quantum groups obtained from semisimple 
complex Lie algebras, the \emph{Quantized Enveloping Algebras}, 
and so we adopt a similar progression of ideas. 

The classical Lie algebra we study, $\sl_2$, is rank one and has 
generators and relations 
\[
[H, X_\pm] = \pm 2 X_\pm, \qquad 
[X_+, X_-] = H 
\]
Let $q$ be a nonzero parameter. We define $U_q(\sl_2)$ as the 
noncommutative algebra generated by 1 and $X_+, X_-, q^{\frac{H}{2}}, 
q^{-\frac{H}{2}}$, with the relations 
$q^{\pm \frac{H}{2}} q^{\mp \frac{H}{2}} = 1$ suggested by the 
notation, together with the nontrivial relations 
\[
q^{\frac{H}{2}} X_{\pm} q^{-\frac{H}{2}} = q^{\pm 1} X_\pm, \qquad 
[X_+, X_-] = \frac{q^H - q^{-H}}{q - q^{-1}} 
\]
The following equations give this algebra the structure of a Hopf algebra: 
\begin{equation} \label{sl2-hopf-1} 
\tri q^{\pm \frac{H}{2}} = q^{\pm \frac{H}{2}} \otimes q^{\pm \frac{H}{2}}, 
\quad 
\tri X_\pm = X_\pm \otimes q^{\frac{H}{2}} + q^{-\frac{H}{2}} \otimes X_\pm 
\end{equation} 
\begin{equation} \label{sl2-hopf-2} 
\eps q^{\pm \frac{H}{2}} = 1, \quad 
\eps X_\pm = 0, \quad 
SX_\pm = -q^{\pm 1} X_\pm, \quad 
Sq^{\pm \frac{H}{2}} = q^{\mp \frac{H}{2}}
\end{equation} 
If we work over the ring $\C[[q]]$ of formal power series in the indeterminate 
$q$, then the Hopf algebra defined by eqns.~\eqref{sl2-hopf-1}-\eqref{sl2-hopf-2}
is quasitriangular, with the following explicit universal R-matrix: 
\[
R = q^{\frac{H \otimes H}{2}} \sum_{n=0}^\infty 
\frac{(1-q^{-2})^n}{[n]!} (q^{\frac{H}{2}} X_+ \otimes q^{-\frac{H}{2}} X_-)^n q^{\frac{n(n-1)}{2}} 
\]
where we have made use of the notation $[n] = \frac{q^n - q^{-n}}{q-q^{-1}}$ 
and $[n]! = [n][n-1]\ldots [1]$. 
This completely specifies the algebra structure and its quantum R-matrix. 
In principle, one at this point should verify all of the Hopf algebra axioms. 
We will not do this here, but refer the reader to the literature \cite{majid}. 

\subsection{Real Forms} 

The classical Lie algebra $\sl_2$ has two inequivalent real forms: 
$\su_2$ and $\su(1,1) \cong \sl(2,\R)$. In this section, we describe the 
corresponding theorem in the case of quantum groups, without giving proofs. 

In the classical case, one can define a real-form of a complex simple 
Lie algebra $\g$ to be an anti-linear anti-involution $\ast$ on $\g$. 
The different real subalgebras which complexify to $\g$ are then identified 
with the $\ast$-invariant subspaces. The $\ast$-structure on $\g$ 
corresponds to a Hopf-$\ast$ structure on the universal enveloping 
algebra $\U(\g)$, so viewing real forms of classical Lie algebras in terms 
of $\ast$ operations abstracts directly to the case of quantum groups. 
Thus we define a real form of a Hopf algebra over $\C$ to be a specification 
of a Hopf $\ast$ structure. In particular, this now gives a well-defined meaning 
to a \emph{unitary representation}, the latter being defined as a 
representation which is $\ast$-equivariant. 

Given these definitions, one finds a direct analogy: just as 
is the case for classical $\sl_2$,  there are two inequivalent real 
forms of the quantum group $U_q(\sl_2)$.

\section{The Quantized Enveloping Algebra $U_q(\g)$.} 

Let $\g$ be a finite dimensional complex 
semisimple Lie algebra $\g$ of rank $N$, 
with Cartan matrix $(a_{ij})$. In this situation, the matrix $a_{ij}$ 
is \emph{symmetrizable}, which means that there are relatively
prime positive integers $d_1, \ldots, d_N$ with the property that the 
matrix $(d_i a_{ij})$ is symmetric. 
Let $\{\alpha_1, \ldots, \alpha_N\}$ be the simple roots, and 
$Q$ the root lattice. Let $(\ ,\ )$ be the positive definite symmetric 
bilinear form on the root space defined by the equations 
\[
(\alpha_i, \alpha_j) = d_i a_{ij} 
\]
which can be realized by dualizing the Killing form $K$. 
We choose \emph{Cartan-Chevalley generators} $H_i, X_{\pm i}$ according
to the following (completely classical) prescription: 
\[
\alpha_i(H_j) = K(d_i H_i, H_j) = a_{ji}, \quad 
K(X_{+i}, X_{-j}) = d_i^{-1} \delta_{ij}
\]
so that we have relations 
\[
[H_i, H_j] = 0, \quad 
[H_i, X_{\pm j}] = \pm a_{ij} X_{\pm j}, \quad 
[X_{+i}, X_{-j}] = \delta_{ij} H_i 
\]
To deform this into the quantized enveloping algebra, we 
can either associate to this root system a set of 
generators $\{q_i^{\pm H_i/2}, X_i, X_{-i} \}$ of a Hopf 
algebra over $\C$ (as we have done for $\sl_2$ previously 
in Section \ref{sec:sl2}) or we can define 
$U_q(\g)$ over the formal power series ring 
$\C[[t]]$ with $q = e^{t/2}, q_i = q^{d_i}$. 

Since the construction involving the $q_i^{\pm H_i/2}$ 
as generators was illustrated for the $\sl_2$ case, 
we will now illustrate the general case using the 
second approach, based on formal power series. 
Accordingly, we define the relations of $U_q(\g)$ as 
\begin{equation} \label{gen-hopf-1} 
[H_i, H_j] = 0, \quad 
[H_i, X_{\pm j}] = \pm a_{ij} X_{\pm j}, \quad 
[X_{+i}, X_{-j}] = \delta_{ij} \frac{q_i^{H_i} - q_i^{-H_i}}{q_i - q_i^{-1}} 
\end{equation} 
\begin{equation} \label{gen-hopf-2} 
\sum_{k=0}^{1-a_{ij}} (-1)^k \qbinom{1-a_{ij}}{k}{q_i} 
X_{\pm i}^{1-a_{ij}-k} X_{\pm j} X_{\pm i}^k = 0, \ 
\text{ for all } i \ne j 
\end{equation} 
where the expressions $\qbinom{n}{m}{q_i}$ are called 
\emph{$q$-binomial coefficients} and has the same 
definition as the usual binomial coefficient, except that 
the factorial function is replaced by the quantum 
factorial $[n]_{q_i}! = [n]_{q_i} [n-1]_{q_i} \dots [1]_{q_i}$ 
where $[n]_{q_i} := \frac{q_i^{n} - q_i^{-n}}{q_i - q_i^{-1}}$. 
\footnote{ quantum factorials $[n]_{q_i}!$ of the form 
utilised here were defined by Heine in 1846!} 
Finally, the coproduct, counit, and antipode are defined by 
the equations 
\begin{equation} \label{gen-hopf-3} 
\tri H_i = H_i \o 1 + 1 \o H_i, \quad 
\tri X_{\pm i} = X_{\pm i} \o q_i^{H_i/2} + q_i^{H_i/2}  \o X_{\pm i},
\end{equation} 
\begin{equation} \label{gen-hopf-4} 
\eps(H_i) = \eps(X_{\pm i}) = 0, \quad SH_i = -H_i, \quad 
SX_{\pm i} = -q_i^{\pm 1} X_{\pm i} 
\end{equation} 

We denote the algebra defined by the generators and relations 
\eqref{gen-hopf-1}-\eqref{gen-hopf-4} by $U_q(\g)$. 
As in the case of the quantized $\sl_2$, we omit the proof of 
these relations and merely state one of its most fundamental properties: 

\begin{proposition} 
The algebra $U_q(\g)$ defined by \eqref{gen-hopf-1}-\eqref{gen-hopf-4} 
is a Hopf algebra. The action of the antipode is $-1$ times
the conjugation by $\rho^\vee$, where $\rho$ is $\frac{1}{2}$ times
the sum of the positive roots, and $\vee$ denotes dualization 
with respect to the Killing form. 
\end{proposition}

\section{Quantum Groups and Braid Groups} 

\begin{definition} 
The \emph{braid group on $n$ strands} ($n \geq 3$) is the group 
$B_n$ generated by $n-1$ elements $\sigma_1, \ldots, \sigma_{n-1}$ 
with the relations 
\[
\sigma_i \sigma_j = \sigma_j \sigma_i \quad \text{ if } \quad \abs{i-j} > 1 
\]
and 
\[
\sigma_i \sigma_{i+1} \sigma_i = \sigma_{i+1} \sigma_i \sigma_{i+1}
\quad \text{ for } \quad 1 \leq i,j \leq n-1
\]
A trivial case is $n = 2$, in which case $B_2 = \Z$ is free on 
one generator. 
\end{definition} 

Proposition \ref{prop:YB-braid} gives a connection between the Yang-Baxter 
equation and linear representations of the braid group.

\begin{proposition} \label{prop:YB-braid} 
Let $V$ be a vector space and $c \in GL(V)$ a linear automorphism 
satisfying the Yang-Baxter equation: 
\[
(c \o \id) (\id \o c) (c \o \id) = (\id \o c) (c \o \id) (\id \o c) 
\]
Define, for $1 \leq i \leq n-1$, linear automorphisms $c_i$ of 
$V^{\o n}$ by: 
\begin{eqnarray*} 
c_1 &=& c \o \id^{\o (n-2)}, \\ 
c_i &=& \id^{\o (i-1)} \o c \o \id^{\o (n-i-1)}, \quad \text{ if } \quad 
1 < i < n-1 \\ 
c_{n-1} &=& \id^{\o (n-2)} \o c 
\end{eqnarray*}
Then there is a unique homomorphism from $B_n$ to the group of linear 
automorphisms of $V^{\o n}$ sending $\sigma_i$ to $c$. 
\end{proposition} 

Thus in any quasitriangular (or ``braided'') Hopf algebra $H$, 
the universal R-matrix allows to define braid group representations
commuting with the adjoint action of $H$
in the tensor powers of any $H$-module. Here, the adjoint action 
of $H$ is an action of $H$ on itself, defined by the equation 
\[
\ad(x)(y) = \sum x_{(1)} y S(x_{(2)}) 
\]
in which the subscripts are a notation which is quite common:
we avoid giving separate names to the homogeneous elements which make
up a general element in a tensor power, and simultaneously take away 
the index labeling the sum, even though there is still a sum. 
Thus, if $R \in V^{\o 2}$ is a non-homogeneous
element in a second tensor power, then the sane notation would be 
$R = \sum_{i=1}^n \alpha_i \otimes \beta_i$. We will follow convention and 
instead write $R = \sum R_{(1)} \otimes R_{(2)}$ with the understanding 
that the symbols $R_{(i)}$ have no independent meaning. 

Braid groups are important in modern representation theory, 
as is illustrated by our discussion of the Drinfeld-Kohno theorem 
(Theorem \ref{thm:Drin-Ko}) in Section \ref{sec:KZeqns}.

\section{Application: the Knizhnik-Zamolodchikov Equations} 
\label{sec:KZeqns} 

Let $\g$ be a complex semisimple Lie algebra as usual, and consider 
$\phi(z) = \phi(z_1, \ldots, z_n$ to be a function of $n$ complex variables
with values in the $n$-fold tensor product of $\g$-modules 
\[
M := M_1 \otimes M_2 \otimes \dots \otimes M_n
\]
Fix an invariant 
symmetric inner product on $\g$, which can be identified with a certain
tensor $\Omega \in \g \otimes \g$. For $i \ne j$, we let 
$\Omega_{ij}$ denote the endomorphism of $M$ induced by letting 
the inner product $\Omega$ act on the $M_i$ and $M_j$ factors, keeping all 
other factors fixed. 

Physicists discovered \cite{KZ} that the correlation functions 
of Wess-Zumino-Witten models (two-dimensional conformal field theory) 
satisfy the following system of differential equations: 
\begin{equation} \label{eqn:KZ} 
\frac{\d\phi(z)}{\d z_i} = \frac{1}{\kappa} 
\sum_{j \ne i} \frac{\Omega_{ij} \phi(z)}{z_i - z_j} 
\end{equation} 
These are the famous Knizhnik-Zamolodchikov equations, after 
which this section is titled. 

The important Kohno-Drinfeld theorem states roughly that the monodromy 
of the solutions of the KZ equations \eqref{eqn:KZ} can be expressed
in terms of the R-matrix of the quantized universal enveloping 
algebra $U_q(\g)$. The quantization parameter $q$ in the quantum group
is connected with the complex constant $\kappa$ which appears 
in the KZ equations by the relation:
\[
q = \exp(2\pi i/\kappa) 
\]

\subsection{The precise statement of the Kohno-Drinfeld theorem}

The following discussion of the Kohno-Drinfeld theorem 
is adapted from \cite{QLBVI}, in which a new (short!) proof of 
the theorem is given, using certain canonical categories 
obtained from the study of Lie bialgebras. 

The category $\Cal O$ for $\g$ is defined 
to be the category of $\h$-diagonalizable 
$\g$-representations, whose weights belong to 
a union of finitely many cones $\lambda-\sum_i \Z_+\alpha_i$, 
$\lambda\in \h^*$, 
and the weight subspaces are finite dimensional. 
Define also the category $\Cal O[[\hbar]]$ of
deformation representations of $\g$, i.e. representations of 
$\g$ on topologically free $k[[h]]$ modules with similar assumptions on the
weights ($\lambda\in \h^*[[\hbar]]$).  

In a similar way one defines the category $\Cal O_\hbar$ 
for the algebra ${\Cal U}$: it is the category 
of ${\Cal U}$-modules which are topologically free over $k[[h]]$ 
and satisfy the same conditions as in the classical case.

\begin{theorem}[Drinfeld-Kohno] \label{thm:Drin-Ko} 
Let $k = \C$. Let 
$V\in {\Cal O}[[\hbar]]$, and 
$V_q=F(V)$ be its image in $\Cal O_\hbar$. Let $F(z_1,...,z_n)$ 
be a function of complex variables $z_1,...,z_n$ 
with values in $ V^{\o n}[\lambda][[\hbar]]$ (the weight subspace of
weight $\lambda$) and consider the system of 
Knizhnik-Zamolodchikov differential equations:
\[
\frac{\d F}{\d z_i}=\frac{\hbar}{2\pi i}\sum_{j\ne
i}\frac{\Omega_{ij}F}{z_i-z_j}.
\]
Then the monodromy representation of the braid group $B_n$ for
this equation is isomorphic to the representation 
of $B_n$ on $V_q^{\otimes n}[\lambda]$ defined by the
formula
\[
b_i\to \sigma_i R_{ii+1},
\]
where $b_i$ are generators of the braid group and 
$\sigma_i$ are the permutation of the i-th and (i+1)-th
components. 
\end{theorem}

\section{Compact Quantum Groups} 

The treatment of Compact Quantum Groups in this paper 
owes much to \cite{woronowicz}. 

\begin{definition} \label{def:CQG} 
Let $G = (A, \Phi)$ where $A$ is a separable unital $C^*$-algebra and 
$\Phi : A \To A \otimes A$ is a unital *-algebra homomorphism. We say that 
$G$ is a \emph{compact quantum group} if 
\begin{enumerate} 
\item The following diagram commutes: 
\[
\xymatrix{ 
A \ar[r]^{\Phi} \ar[d]_{\Phi}  & A \otimes A \ar[d]^{\Phi \otimes \id}  \\ 
A \otimes A \ar[r]^{\id \otimes \Phi} & A \otimes A \otimes A } 
\]
\item  The sets $\{ (b \otimes I)\Phi(c) : b, c \in A\}$ and 
$\{ (I \otimes b)\Phi(c) : b, c \in A\}$ are linearly dense subsets of $A \otimes A$. 
\end{enumerate} 
\end{definition} 

Compact quantum matrix groups can be recovered from the above formalism
in the following way. 
If $(A, u)$ is a compact quantum matrix 
group, and $\Phi$ is the corresponding comultiplication, then $(A, \Phi)$ 
satisfies the axioms given above. 

As one might expect, the simplest case to understand is the case in which 
the algebra $A$ defining the quantum group is commutative. 
In this situation, it follows from Gelfand-Naimark theory that
$\exists$ some compact space $\Lambda$ such that  
$A = C(\Lambda)$. The comultiplication map $\Phi$ can also be understood
within the context of Gelfand-Naimark theory. In this case the implication
is that there exists a continuous map $\Lambda \times \Lambda \To \Lambda$
denoted by $(\lambda_1, \lambda_2) \Mapsto \lambda_1 \cdot \lambda_2$ 
such that for any $a \in C(\Lambda)$ and for any $\lambda_1, \lambda_2$, 
we have 
\[
(\Phi a)(\lambda_1, \lambda_2) = a(\lambda_1 \cdot \lambda_2) 
\]
Moreover, condition (1) of Definition \ref{def:CQG} means that 
the ``multiplication'' $\cdot$ is associative, and the two
density conditions in (2) of Def.~\ref{def:CQG} imply the two-sided
cancellation law for ordinary group multiplication. It follows that 
$\Lambda$ in this situation is a topological group. 
This justifies the notion that a compact quantum group is a 
generalization of the classical idea of a compact group to an 
underlying ``noncommutative space,'' and as such, fits into Alain 
Connes' formulation of non-commutative geometry.

\section{Haar Measure on Quantum Groups} 

Let $G = (A, \Phi)$ be a compact quantum group. The goal of  this
section will be to describe the generalization of Haar measure to this
case. First it is necessary to define the convolution product. 
In the following definitions, $\xi, \xi'$ are continuous linear
functionals on $A$, and $a \in A$. One can convolve an element of $A$ 
with a functional to produce an element of $A$, and one 
can convolve two functionals to produce another functional.  
The precise statements are: 
\begin{eqnarray*} 
\xi \ast a &:=& (\id \otimes \xi) \Phi(a)  \\ 
a \ast \xi &:=& (\xi \otimes \id) \Phi(a) \\ 
\xi' \ast \xi &:=& (\xi' \otimes \xi)\Phi 
\end{eqnarray*} 

By definition, a \emph{state} $\xi$ on a $C^*$-algebra $A$ 
is a positive normalized linear functional on $A$.
This terminology comes from the relationship which exists 
between quantum mechanics and $C^*$-algebras, provided 
by the Gelfand-Naimark-Segal (GNS) construction. Briefly,
from the data $(A, \xi)$, the
GNS construction gives a Hilbert space $\H$, a unit cyclic 
vector $\psi \in \H$, and a $\ast$-representation 
$\pi_\xi : A \to \B(\H)$ such that 
\[
\xi(a) = \left< \pi_\xi(a)\psi, \psi \right> 
\text{ for all } a \in A 
\]

The following fundamental result establishes the analogue 
of Haar measure for compact quantum groups. 

\begin{theorem} 
Let $G = (A, \Phi)$ be a compact quantum group. Then there exists 
a unique state $h$ on $A$ satisfying 
\begin{equation} \label{eqn:Haar} 
a \ast h = h \ast a = h(a) I \text{ for all } a \in A
\end{equation} 
\end{theorem} 

To relate this to the classical theory of Haar measure on a compact 
group, one must notice that if $A = C(\Lambda)$ is the algebra of 
continuous, complex-valued functions on a compact 
topological group $\Lambda$,
then the Haar measure $\mu$ on $\Lambda$ is certainly a 
positive measure 
on the Borel $\sigma$-field of measurable sets in $\Lambda$. However, 
it can be equivalently described as a functional on $A$ 
given by $\mu(f) := \int_\Lambda f \, d\mu$. 
Convolutions of measures with test functions can be defined as in 
distribution theory, and it is straightforward to show that 
\eqref{eqn:Haar} holds for the Haar measure (state) $\mu$. 

\section{Representations of Compact Quantum Groups} 

Before we can continue, it is necessary to 
introduce some preliminary notation and terminology from 
operator algebra theory. Let $\A$ be a $C^*$-algebra and let $A, B$ be 
bounded linear operators on $\A$. We say that $B$ is the adjoint of $A$, 
if $A(x)^* \, y = x^* \, B(y)$ for all $x,y \in \A$. In this situation 
we denote $B$ by $A^*$. The \emph{multiplier algebra} $M(\A)$ 
is defined to be the subalgebra of the bounded operators $B(\A)$, 
defined by 
\[
M(\A) = \{ A \in B(\A) \mid A^* \text{ exists } \} 
\]

\begin{definition}
Let $G = (A, \Phi)$ be a compact quantum group and $\H$ a Hilbert space.
Let $K$ denote the algebra of bounded, compact operators on $\H$. 
A \emph{(strongly continuous) unitary representation} of $G$ on $\H$ 
is a  unitary element $v \in M(K \otimes A)$ such that 
\[
(\id \otimes \Phi) v = v_{12} v_{13} 
\]
\end{definition} 

One should think of $v$ ``acting on'' $\H$ in the following sense. 
As an element of the multiplier algebra $M(K \otimes A)$, 
$v$ is a bounded operator on $K \otimes A$ possessing an adjoint 
with respect to the $\ast$-product on $K \otimes A$ induced by the 
tensor product construction. 
Thus, using $v$, every element $a \in A$ ``acts on'' the compact 
operators $K$ via the pairing $(a, k) \To v(k \otimes a)$. 

If the Hilbert space $\H$ is just $\C^N$, then the multiplier 
algebra of $K \otimes A$ is $M_N(A)$. 
A finite-dimensional representation of a compact quantum group is 
therefore by definition, a unitary matrix (with matrix elements from
the algebra $A$)  
\[
v = (v_{k\ell})_{k,\ell = 1,\ldots, N}
\]
satisfying the condition that 
\begin{equation} \label{eqn:fin-dim-rep}
\Phi(v_{k\ell}) = \sum_r v_{kr} v_{r\ell} \text{ for all } 
k, \ell = 1, \ldots, N 
\end{equation} 
In other words, matrix multiplication is the quantum ``group law'' 
on $G$, determined by the comultiplication. 
Equation \eqref{eqn:fin-dim-rep} is the analogue 
(with appropriate arrows reversed) of the classical 
notion that for a matrix representation of a group $\rho : G \to GL(V)$, 
the matrix elements of products are given by the formula 
$\rho(g\,h)_{k\ell} = \sum_r \rho(g)_{kr} \,\rho(h)_{r\ell}$.

It is a remarkable result that even at this level of generality, 
one obtains a complete reducibility result for these representations. 

\begin{theorem} 
Let $v$ be a unitary representation of a compact quantum group 
$G = (A, \Phi)$ on any Hilbert space $\H$. 
Then $v$ decomposes as a direct sum of finite-dimensional 
irreducible representations. 
\end{theorem}

\end{document}